\documentclass{article}

\RequirePackage{hyperref}
\usepackage{amsmath}
\usepackage{amssymb}
\usepackage{tikz}
\usepackage{listings}
\usepackage{tikz}
\usetikzlibrary{positioning,shapes,arrows,automata,decorations.pathreplacing,angles,quotes}

\usepackage{subcaption}
\usepackage{caption}
\usepackage{wrapfig}

\newcommand{\R}{\mathbb{R}}
\newcommand{\V}{{\bf v}}
\newcommand{\X}{{\bf x}}

\newcommand{\Y}{{\bf y}}
\renewcommand{\P}{{\bf p}}

\newcommand{\ass}{{\scriptsize \, \mathrel{{:}{=}}\,}}
\newcommand{\fma}{{\small {\tt fma}}}

\newcommand{\pe}{{\scriptsize \, \mathrel{{+}{=}}\, }}

\newcommand\Ert{{\footnotesize \textsc{Cost}}}

\newtheorem{Rule}{Rule}
\newtheorem{Problem}{Problem}
\newtheorem{example}{Example}
\newtheorem{Definition}{Definition}

\begin{document}

\title{Elimination Techniques for Algorithmic Differentiation Revisited}

\author{Uwe Naumann\thanks{corresponding author: {\tt naumann@stce.rwth-aachen.de}}
\and Erik Schneidereit \and Simon M\"artens \and Markus Towara\thanks{all: Informatik 12: Software and Tools for Computational Engineering, RWTH Aachen University, 52056 Aachen, Germany.}}

\date{}

\maketitle

\begin{abstract}
All known elimination techniques for (first-order) algorithmic differentiation 
(AD) rely on Jacobians to be given for a set of relevant elemental functions. 
Realistically, elemental tangents and adjoints are given instead. They can
be obtained 
by applying software tools for AD to the parts of a given modular numerical 
simulation.
The novel generalized face elimination rule
proposed in this article 
facilitates the rigorous exploitation of 
associativity of the chain rule of differentiation
at arbitrary levels of granularity ranging from 
elemental scalar (state of the art) to multivariate 
vector functions with given elemental tangents and adjoints. 
The implied combinatorial {\sc Generalized Face Elimination} 
problem asks for a face elimination sequence of minimal computational cost. 
Simple branch and bound and greedy heuristic methods are employed as
a baseline for further research into more powerful algorithms motivated by 
promising first test results. The latter can be reproduced with 
the help of an open-source reference implementation.
\end{abstract}

\section{Introduction}

Research in elimination methods for algorithmic differentiation (AD) \cite{Griewank2008EDP,Naumann2012TAo}
has been targeting decompositions of differentiable numerical programs
into simple scalar elemental functions (mainly the built-in numerical 
operators and intrinsic functions provided by modern programming languages) 
for more than three decades. Elimination methods aim to exploit associativity
of the chain rule of differentiation in order to minimize the computational cost of accumulating Jacobians, Hessians, or possibly even higher derivative tensors.
Applicability
is typically restricted to static code fragments\footnote{Data dependences are defined uniquely at compile time.} of limited complexity
(e.g. basic blocks). Theoretical benefits due to lower operations counts
require careful implementation in order to translate into actual run time
savings. The application of available optimization methods for elimination
methods to real-world large-scale AD 
missions is hindered substantially by the present restrictions. 

This paper aims to address these issues by lowering the level of
granularity to multivariate vector elementals of potentially significant
computational cost. 
Additionally, all known elimination methods are generalized to
the practically relevant scenario, where elemental tangent and adjoint models
are available instead of elemental Jacobians. Illustration is provided by the 
following example to be referred to repeatedly throughout this article.
\begin{example} \label{ex:Newton1}
Newton's method for the iterative approximation of a solution for 
parameterized systems of nonlinear equations $R(\X,\P)=0$
with twice 
continuously differentiable residual 
$R : \R^n \times \R^m \rightarrow \R^n$
evolves for fixed $\P$ as 
$\X^{i+1}=\X^i-R'(\X^i,\P)^{-1} \cdot R(\X^i,\P)$ 
starting from a given $\X^0$ and where $R'$ denotes the invertible Jacobian 
of $R$ with respect to the state $\X.$ Termination after $k$ steps may be due 
to the 
requirement to reach a given accuracy, e.g. $|R(\X^k,\P)|\leq \epsilon,$ or
$k$ may have to be fixed in order to meet constraints on the overall 
computational cost.

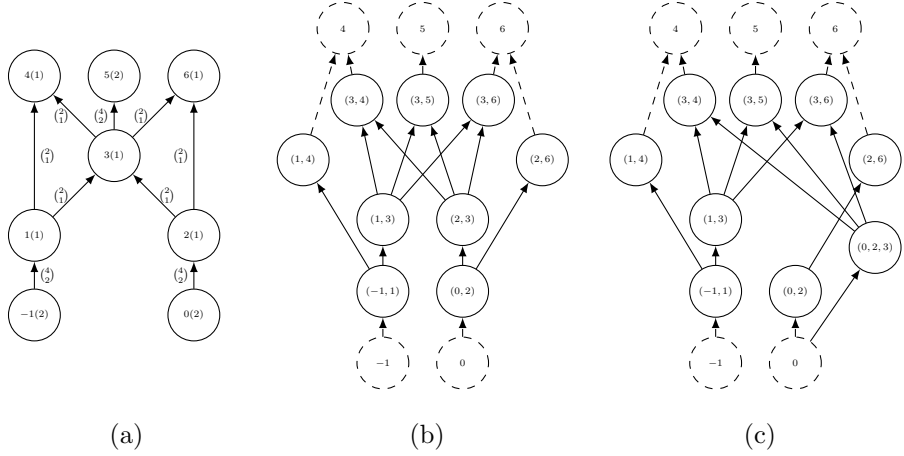
\begin{figure*}
	\begin{tabular}{ccc}
\begin{minipage}[c]{.26\linewidth}
\begin{tikzpicture}[scale=.53, transform shape]
  \begin{pgfscope}
    \tikzstyle{every node}=[draw,circle,minimum size=1.3cm]
	  \node (-1) at (0,0) {\scriptsize $-1 (2)$};
  \node (0) at (4,0) {\scriptsize $0 (2)$};
  \node (1) at (0,2) {\scriptsize $1 (1)$};
  \node (2) at (4,2) {\scriptsize $2 (1)$};
	  \node (3) at (2,4) {\scriptsize $3 (1)$};
    \node (4) at (0,6) {\scriptsize $4 (1)$};
	  \node (5) at (2,6) {\scriptsize $5 (2)$};
    \node (6) at (4,6) {\scriptsize $6 (1)$};
  \end{pgfscope}
 \begin{scope}[-latex]
	 \draw (-1) -- (1) node[midway,right] {\footnotesize $\binom{4}{2}$}; 
	 \draw (0) -- (2) node[midway,left] {\color{black} \footnotesize $\binom{4}{2}$};
 \draw (1) -- (3) node[midway,left] {\footnotesize $\binom{2}{1}$};
 \draw (1) -- (4) node[midway,right] {\footnotesize $\binom{2}{1}$};
	 \draw (2) -- (3) node[midway,right] {\color{black} \footnotesize $\binom{2}{1}$};
	 \draw (2) -- (6) node[midway,left] {\footnotesize $\binom{2}{1}$};
 \draw (3) -- (4) node[midway,left] {\footnotesize $\binom{2}{1}$};
 \draw (3) -- (5) node[midway,left] {\footnotesize $\binom{4}{2}$};
	 \draw (3) -- (6) node[midway,left] {\footnotesize $\binom{2}{1}$};
  \end{scope}
\end{tikzpicture} 
\end{minipage}
	&
\begin{minipage}[c]{.33\linewidth}
\begin{tikzpicture}[scale=.53, transform shape]
  \begin{pgfscope}
    \tikzstyle{every node}=[draw,circle,minimum size=1.3cm]
	  \node[dashed] (-1d) at (2,0.4) {\scriptsize $-1$};
	  \node[dashed] (0d) at (4,0.4) {\scriptsize $0$};
	  \node (-11) at (2,2.2) {\scriptsize $(-1,1)$};
	  \node (02) at (4,2.2) {\scriptsize $(0,2)$};
	  \node (13) at (2,4) {\scriptsize $(1,3)$};
	  \node (23) at (4,4) {\scriptsize $(2,3)$};
	  \node (14) at (0,5.5) {\scriptsize $(1,4)$};
	  \node (26) at (6,5.5) {\scriptsize $(2,6)$};
	  \node (34) at (1.35,7) {\scriptsize $(3,4)$};
	  \node (35) at (3,7) {\scriptsize $(3,5)$};
	  \node (36) at (4.65,7) {\scriptsize $(3,6)$};
	  \node [dashed](4d) at (1,8.8) {\scriptsize $4$};
	  \node [dashed](5d) at (3,8.8) {\scriptsize $5$};
	  \node[dashed] (6d) at (5,8.8) {\scriptsize $6$};
  \end{pgfscope}
 \begin{scope}[-latex]
	 \draw [dashed](-1d) -- (-11);
	 \draw [dashed](0d) -- (02);
	 \draw (-11) -- (13);
	 \draw (-11) -- (14);
	 \draw (02) -- (23);
	 \draw (02) -- (26);
	 \draw (13) -- (34);
	 \draw (13) -- (35);
	 \draw (13) -- (36);
	 \draw (23) -- (34);
	 \draw (23) -- (35);
	 \draw (23) -- (36);
	 \draw[dashed] (14) -- (4d);
	 \draw[dashed] (34) -- (4d);
	 \draw[dashed] (35) -- (5d);
	 \draw[dashed] (36) -- (6d);
	 \draw[dashed] (26) -- (6d);
  \end{scope}
\end{tikzpicture} 
\end{minipage} 
&
\begin{minipage}[c]{.33\linewidth}
\begin{tikzpicture}[scale=.53, transform shape]
  \begin{pgfscope}
    \tikzstyle{every node}=[draw,circle,minimum size=1.3cm]
	  \node[dashed] (-1d) at (2,0.4) {\scriptsize $-1$};
	  \node[dashed] (0d) at (4,0.4) {\scriptsize $0$};
	  \node (-11) at (2,2.2) {\scriptsize $(-1,1)$};
	  \node (02) at (4,2.2) {\scriptsize $(0,2)$};
	  \node (13) at (2,4) {\scriptsize $(1,3)$};
	  \node (023) at (6,3.3) {\scriptsize $(0,2,3)$};
	  \node (14) at (0,5.5) {\scriptsize $(1,4)$};
	  \node (26) at (6,5.5) {\scriptsize $(2,6)$};
	  \node (34) at (1.35,7) {\scriptsize $(3,4)$};
	  \node (35) at (3,7) {\scriptsize $(3,5)$};
	  \node (36) at (4.65,7) {\scriptsize $(3,6)$};
	  \node [dashed](4d) at (1,8.8) {\scriptsize $4$};
	  \node [dashed](5d) at (3,8.8) {\scriptsize $5$};
	  \node[dashed] (6d) at (5,8.8) {\scriptsize $6$};
  \end{pgfscope}
 \begin{scope}[-latex]
	 \draw [dashed](-1d) -- (-11);
	 \draw [dashed](0d) -- (02);
	 \draw (-11) -- (13);
	 \draw (-11) -- (14);
	 \draw (0d) -- (023);
	 \draw (023) -- (34);
	 \draw (023) -- (35);
	 \draw (023) -- (36);
	 \draw (02) -- (26);
	 \draw (13) -- (34);
	 \draw (13) -- (35);
	 \draw (13) -- (36);
	 \draw[dashed] (14) -- (4d);
	 \draw[dashed] (34) -- (4d);
	 \draw[dashed] (35) -- (5d);
	 \draw[dashed] (36) -- (6d);
	 \draw[dashed] (26) -- (6d);
  \end{scope}
\end{tikzpicture} 
\end{minipage} \\
\\
(a) & (b) & (c)
\end{tabular}
	\caption{(a): $G$; \hspace{.2cm} (b): $\tilde{G}$; \hspace{.2cm} (c) $\tilde G-(0,2,3)$} \label{fig:bat}
\end{figure*}
Think of an inexact Newton method \cite{Dembo1982INM}, where the
three elemental functions
$F_1 : \R^n \times \R^m \rightarrow \R^{n \times n},$
$F_2 : \R^n \times \R^m \rightarrow \R^n,$
$F_3 : \R^n \times \R^{n \times n} \times \R^n \rightarrow \R^n$
represent the approximation of the Jacobian of the residual
($F_1=R'$) by finite differences, the evaluation of the residual ($F_2=R$) 
and the iterative approximation of the solution of
the linear Newton system including incrementation of the state
($F_3=\X^i+F_1^{-1} \cdot F_2$). An implicit function $\X=\X(\P)$ is defined
by $R(\X,\P)=0.$ We are interested in tangents and/or adjoints of $\X$ with 
respect to $\P.$\footnote{Without loss of generality, elimination techniques are
discussed in the context of the accumulation of full Jacobians. The generalization to
tangents and adjoints beyond products of the Jacobian with identities follows
seamlessly.} Both can (and should) be derived symbolically at the (exact) 
solution \cite{Gilbert1992Ada}.
AD is to be preferred in the inexact case as this 
method provides exact (with machine accuracy) tangents and adjoints of the 
actually computed $\X^k$. 
\end{example}
AD mission planning refers to the task of determining for a given decomposition
of a modular numerical simulation program into elemental subprograms an order 
of evaluations of the elemental tangent and adjoint subprograms such that
the desired derivative is obtained at minimal computational cost.
It relies on realistic cost estimates for all elemental tangents and 
adjoints. The resulting optimized derivative code amounts to the evaluation of
a generalized face elimination sequence.
E.g., in \cite{Forth2004JCG} it was 
shown that theoretical savings in computational cost
due to the exploitation of associativity of the chain rule 
translate into corresponding improvements in terms of run time. 
Similar observations were made in \cite{Naumann2022HCB} in the context of 
computing Hessians.

\section{State of the Art} \label{sec:2}
We consider implementations of multivariate vector functions
\begin{equation} \label{eqn:F}
        F: \R^n \rightarrow \R^m : \X \mapsto \Y \ass F(\X)
\end{equation}
over the real (floating-point) numbers $\R$
as differentiable computer programs with Jacobian matrices (also: Jacobians)
\begin{equation} \label{eqn:derivs}
        F'=F'(\X) \equiv \frac{\partial F}{\partial \X}(\X) \in \R^{m \times n} \; .
\end{equation}
{\em Figure~\ref{fig:bat}} will provide graphical illustration for the concepts to be discussed in the following.
We use lower case boldface letters to distinguish vectors from (non-bold)
scalars. Matrices as well as multivariate vector functions are denoted by upper case 
letters. 
AD yields {\em scalar tangents} 
\begin{equation} \label{eqn:tan}
\R^m \ni \dot{\Y} \ass F'(\X) \cdot \dot{\X} 
\end{equation} 
and {\em scalar adjoints} 
\begin{equation} \label{eqn:adj}
\R^{1 \times n} \ni \bar{\X} \ass \bar{\Y} \cdot F'(\X) \; .
\end{equation} 
The full Jacobian can be obtained (also: harvested) by setting 
$\dot{\X}$ or $\bar{\Y}^T$ equal to (also: seeding) the Cartesian basis
vectors in $\R^n$ or $\R^m,$ respectively, implying cheap gradients ($m=1$)
\cite{Griewank2012Wit}.
Reapplication results in second- and higher-order tangents and adjoints. 
Derivatives of $F$ 
can be evaluated efficiently and with machine accuracy.
We use $=$ to denote mathematical equality, $\equiv$ in the sense of ``is defined as'' and $\ass$ to represent assignment according to imperative programming.

$F$ is referred to as the {\em primal program}. It evaluates a
partially
ordered
sequence of continuously differentiable {\em elemental functions} (also: {\em elementals})
$$F_j=F_j(\V_i)_{i \prec j} : \R^{\sum_{i \prec j} n_i} \rightarrow \R^{n_j}$$
as a {\em single assignment code}
\begin{equation} \label{eqn:sac}
\R^{n_j} \ni \V_j \ass F_j(\V_i)_{i \prec j} \quad
\text{for}~j=1,\ldots,q
\end{equation}
and where, adopting the notation from \cite{Griewank2008EDP}, $i \prec j$ if and only if $\V_i$ is an argument of $F_j.$
Hence, $\X=(\V_i)_{i\in X}$ and $\Y=(\V_i)_{i \in Y},$
where 
$X=\{1-n,\ldots,0\}$ and $Y =\{p+1,\ldots,q\}$ for
$p \equiv q-m.$
In the following we distinguish between scalar ($n_j=1$ for $j=1,\ldots,q$)
and vector ($n_j \geq 1$) scenarios.
Without loss of generality (w.l.o.g.), we assume the results $\Y$ of $F$ to be
represented by the last $m$ entries of $\V.$ 

Elementals represent differentiable subprograms of
the primal program with given {\em elemental Jacobians}
\begin{equation} \label{eqn:eljac}
F'(i,j) \equiv \frac{\partial F_j}{\partial \V_i}(\V_k)_{k \prec j} \in \R^{n_j \times n_i} \; .
\end{equation}
Augmentation of Equation~(\ref{eqn:sac}) with the computation of all elemental Jacobians as in Equation~(\ref{eqn:eljac}) yields the {\em augmented single assignment code}.

\noindent The state of the art in elimination methods for AD has been targeting the 
scalar scenario.
The corresponding augmented single assignment code induces an edge-labeled directed acyclic graph (dag) $G=(V,E)$
with integer vertices $V=\{1-n,\ldots,q\}$ representing the (scalar) $v_j$ and edges $E \subseteq V \times V =\{(i,j) : i \prec j\}.$
Minimal (without predecessors) vertices representing the scalar entries of $\X$ are 
collected in $X \subseteq V.$ The maximal vertices (without successors) 
$Y \subseteq V$ represent the scalar entries of $\Y.$ 
Scalar first derivatives $F'(i,j)$ are associated with all edges $(i,j) \in E.$ 

\begin{example} \label{ex:1}
The dag of the primal C-program 
\begin{verbatim}
void F(float *x, float *y) {
  a=sin(x[0]); b=cos(x[1]); c=a/b; 
  y[0]=a*c; y[1]=exp(c); y[2]=c+b;
}
\end{verbatim}
is depicted in {\em Figure~\ref{fig:bat}~(a).} 
\end{example}
The chain rule of differentiation becomes
\begin{equation} \label{eqn:cr}
\frac{d v_t}{d v_s} = \sum_{\pi=(s,*,t)} \prod_{(i,j) \in \pi} F'(i,j)
\end{equation}
\cite{Baur1983TCo}
where the wildcard $*$ represents the possibly empty remainder of (vertex) 
paths connecting $s \in X$ and $t \in Y$ in $G;$ e.g. 
$(2,3) \in *$ for $(0,*,6)$ in Figure~\ref{fig:bat}~(a). 
Associativity of the chain rule implies that Equation~(\ref{eqn:cr}) holds for any 
subgraph of $G.$ 

Vector dags follow naturally for $n_j \geq 1.$ The $F'(i,j) \in \R^{n_j \times n_i}$ become elemental Jacobian matrices. 
A vector version of the dag in Figure~\ref{fig:bat}~(a) is used for illustration
throughout
this article. Vertices are hence annotated with vector sizes (in parentheses).
Edges $(i,j)$ are annotated with elapsed run times 
for evaluating a scalar tangent 
($\Ert(\dot{F}(i,j),$ bottom number) and a scalar adjoint 
($\Ert(\bar{F}(i,j),$ top number); see Section~\ref{sec:3} for further details.
\begin{figure}
\centering
\begin{tikzpicture}[scale=.53, transform shape]
\begin{pgfscope}
  \tikzstyle{every node}=[draw,circle,minimum size=1.3cm]
  \node (-1) at (0,0) {\scriptsize $-1$};
  \node (0) at (4,0) {\scriptsize $0$};
  \node (1) at (2,3) {\scriptsize $1$};
  \node (2) at (5,3.5) {\scriptsize $2$};
  \node (3) at (0,6) {\scriptsize $3$};
\end{pgfscope}
\begin{scope}[-latex]
 \draw (-1) -- (1) node[midway,left] {\footnotesize $\binom{4000}{2000}$};
 \draw (0) -- (1) node[midway,right] {\footnotesize $\binom{4000}{2000}$};
 \draw (-1) -- (2) node[midway,left] {\footnotesize $\binom{400}{200}$};
 \draw (0) -- (2) node[midway,right] {\footnotesize $\binom{400}{200}$};
 \draw (-1) -- (3) node[midway,left] {\footnotesize $\binom{0}{0}$};
 \draw (1) -- (3) node[midway,left] {\footnotesize $\binom{2000}{1000}$};
 \draw (2) -- (3) node[midway,above] {\footnotesize $\binom{1600}{800}$};
\end{scope}
\end{tikzpicture} 
\caption{Annotated dag of a single Newton step} \label{fig:newton}
\end{figure}
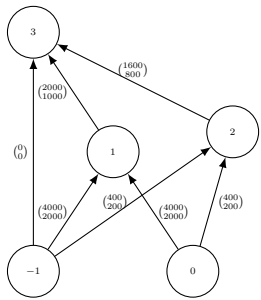
\begin{example} \label{ex:2}
Building on {\sc Example~\ref{ex:Newton1}}, 
the dag of a single Newton step is shown in 
{\em Figure~\ref{fig:newton}}
for
$\V_{-1}=\X^i$ and $\V_0=\P.$
Elemental Jacobians
\begin{align*}
        F'(-1,1)&\equiv\frac{\partial R'}{\partial \X^i} \in \R^{n^2 \times n}, \quad
F'(-1,2)\equiv\frac{\partial R}{\partial \X^i} \in \R^{n \times n} \\
        F'(0,1)&\equiv\frac{\partial R'}{\partial \P} \in \R^{n^2 \times m}, \quad
F'(0,2)\equiv\frac{\partial R}{\partial \P} \in \R^{n \times m} \\
        F'(-1,3)&\equiv\frac{\partial \X^{i+1}}{\partial \X^i}=I_n \in \R^{n \times n} \\
        F'(1,3)&\equiv\frac{\partial \X^{i+1}}{\partial R'} \in \R^{n \times n^2}, \quad 
        F'(2,3)\equiv\frac{\partial \X^{i+1}}{\partial R} \in \R^{n \times n}
        \end{align*}
are associated with the respective edges. 
Equation~(\ref{eqn:cr}) yields
the Jacobian $\frac{d \X^{i+1}}{d \binom{\X^i}{\P}}$ as
$$
\begin{pmatrix}
        F'(-1,3)+F'(1,3) \cdot F'(-1,1)+F'(2,3) \cdot F'(-1,2) \\
        F'(1,3) \cdot F'(0,1)+F'(2,3) \cdot F'(0,2) 
\end{pmatrix}.
$$
\end{example} 
(Cost-)Optimal interpretations of the chain rule along 
single paths amount to the solution of the NP-complete {\sc Matrix Chain 
Product} 
problem \cite{Naumann2020CSC}. The proof presented therein exploits potential
algebraic dependences (more specifically, equality) between the elemental
partial derivatives. Even in the {\em structural} scenario, where partials are assumed to
be algebraically independent, the scope for optimization can be quadratic in 
the problem size. Simply compare the operation counts induced by the 
best (${\mathcal O}(n\cdot q)$)
versus the worst (${\mathcal O}(n^3 \cdot q)$) bracketings for
$F=F_q \circ F_{q-1} \circ \ldots \circ F_2 \circ F_1$ with
$F_{2 \cdot i-1} : \R \rightarrow \R^n$ and
$F_{2 \cdot i} : \R^n \rightarrow \R$ 
while $i=1,\ldots,q/2.$ 
Optimal bracketings for Jacobian chain products can be obtained efficiently 
by dynamic programming \cite{Godbole1973}.


We aim to restrict the amount of material reproduced from \cite{Naumann2004Oao} 
to a minimum. Nevertheless, some essential results about state-of-the-art 
elimination techniques need to be restated explicitly in order to make this 
article as self-contained as possible.

{\em Face elimination} follows immediately from the application of 
Equation~(\ref{eqn:cr}) to paths $(i,j,k),$ $i,j,k \in V,$ of (vertex-)length three. The resulting
modifications in data dependence between intermediate variables cannot
be represented structurally in $G$ in general. For example, 
elimination of $(1,3,4)$ 
in Figure~\ref{fig:bat}~(a) 
must not result in the removal of $(1,3)$ (equivalent to setting $F'(1,3) \ass 0$) as $(1,3) \in (1,3,5)$ and $(1,3) \in (1,3,6).$ Similarly,
$(3,4) \in (2,3,4)$ must remain in $G.$ To resolve this conflict, face 
elimination
is defined on the following variant of the line graph of $G.$
\begin{Definition}[\cite{Naumann2004Oao}] \label{def:dcg}
The {\em line dag} 
$\tilde{G}=(\tilde{V},\tilde{E})$ of $G$ consists of
vertices $\tilde{V}=\tilde{X} \cup \tilde{Z} \cup \tilde{Y}$, such that
$\tilde{X} \cap \tilde{Z} =
\tilde{X} \cap \tilde{Y} =
\tilde{Y} \cap \tilde{Z}
= \varnothing,$ and edges
$\tilde{E}=\tilde{E}_X \cup \tilde{E}_Z \cup \tilde{E}_Y,$ with
$\tilde{E}_X \cap \tilde{E}_Z  =
\tilde{E}_X \cap \tilde{E}_Y =
\tilde{E}_Y \cap \tilde{E}_Z
= \varnothing.$
It is defined by the following construction:
\begin{itemize}
\item[] $\tilde{X}=X; \; \; \tilde{Z}=E; \;\; \tilde{Y}=Y$ 
\item[] $\forall~(i,j),(j,k) \in E~:~((i,j),(j,k)) \in \tilde{E}_Z;$
\item[] $ (i,(i,j)) \in \tilde{E}_X $ for $i \in X$ and $(i,j)$ in $E;$
\item[] $ ((i,j),j) \in \tilde{E}_Y $ for $(i,j) \in E$ and $j \in Y.$
\end{itemize}
\end{Definition}
The line dag for {\sc Example~\ref{ex:1}} is depicted in Figure~\ref{fig:bat}~(b).
Initially, it consists of vertices labeled with index pairs exclusively.
This is due to all vertices in $\tilde{Z}$ representing edges (paths of vertex length 
two) from $G.$
New vertices (fill) are added according to the chain rule to represent
paths of (vertex) length three and more. The corresponding
labels become triplets, quadruples etc. 
For example,
for the dag in {\sc Example~\ref{ex:1}} 
$F'(-1,1,3,6) \equiv F'(3,6) \cdot F'(1,3) \cdot F'(-1,1)$ is computed.
In the following, we write 
$(i,*,j)$ to denote a vertex due to the chain rule applied to some
path connecting a vertex $(i,i')$ with a vertex $(j',j).$ The wildcard
$*$ represents elements from the set of all feasible sub-paths; e.g., $1,3$ for
$F'(-1,*,6)$ as $(i,i')=(-1,1)$ and $(j',j)=(3,6).$ 
The set of predecessors [successors]
of $(i,*,j) \in \tilde{V}$ is denoted as $P(i,*,j) \subset \tilde{V}$ [$S(i,*,j) \subset \tilde{V}$].
We write $(i,*)$ [$(*,j)$] if explicit specification of the
sink [source] is not required in the current context.
\begin{Rule}[\cite{Naumann2004Oao}] \label{r1}
	{\em Face elimination} is defined for all $((i,*,j),(j,*,k)) \in \tilde{E}_Z$ as follows:
\begin{enumerate}
	\item If $\nexists (i,*,k) \in \tilde{V}:$ 
\begin{itemize}
\item[] $P(i,*,k)=P(i,*,j) \, \text{and} \; S(i,*,k)=S(j,*,k)$ 
\end{itemize}
then {\em (fill)}
\begin{itemize}
\item[]	$\tilde{V} \ass \tilde{V} \cup (i,*,j,*,k)$ 
\item[] $\tilde{E} \ass \tilde{E}$ 
\item[] \quad $\cup\; (P(i,*,j) \times (i,*,j,*,k))$ {\em (up-link)}
\item[] \quad $\cup \;((i,*,j,*,k) \times S(j,*,k))$ {\em (down-link)}
\item[] $F'(i,*,j,*,k) \ass 0.$ 
\end{itemize}
\item $F'(i,*,j,*,k) \pe F'(j,*,k) \cdot F'(i,*,j)$ {\em (absorb)}
\item $\tilde{E} \ass \tilde{E} \setminus ((i,*,j),(j,*,k))$ {\em (remove edge)}
\item If $S(i,*,j)=\varnothing$ then {\em (remove isolated vertex)}
\begin{itemize}
	\item[] $\tilde{V} \ass \tilde{V} \setminus (i,*,j)$ 
	\item[] $\tilde{E} \ass \tilde{E} \setminus (P(i,*,j) \times (i,*,j))$ 
\end{itemize}
else if $\exists(i,*) \neq (i,*,j) \in \tilde{V}:~S(i,*)=S(i,*,j)$ 
		then {\em (merge)}
\begin{itemize}
	\item[] $F'(i,*,j) \pe F'(i,*)$ 
	\item[] $\tilde{V} \ass \tilde{V} \setminus (i,*)$ 
	\item[] $\tilde{E} \ass \tilde{E} \setminus ((i,*) \times S(i,*) \cup P(i,*) \times (i,*)).$ 
\end{itemize}
\item If $P(j,*,k)=\varnothing$ then {\em (remove isolated vertex)}
\begin{itemize}
	\item[] $\tilde{V} \ass \tilde{V} \setminus (j,*,k)$
	\item[] $\tilde{E} \ass \tilde{E} \setminus ((j,*,k) \times S(j,*,k))$ 
\end{itemize}
else if $\exists(*,k) \neq (j,*,k) \in \tilde{V}:~P(*,k)=P(j,*,k)$ then {\em (merge)}
\begin{itemize}
	\item[] $F'(j,*,k) \pe F'(*,k)$ 
	\item[] $\tilde{V} \ass \tilde{V} \setminus (*,k).$ 
	\item[] $\tilde{E} \ass \tilde{E} \setminus ((*,k) \times S(*,k) \cup P(*,k) \times (*,k)).$ 
\end{itemize}
\end{enumerate}
\end{Rule}
Fill-in is generated with vanishing derivative $F'(i,*,j,*,k)$ (see {\em 1}) 
followed
by the fused multiply-add (\fma) operation due to the chain rule (see {\em 2}).
The edge $((i,*,j),(j,*,k))$ is removed from $\tilde{G}$ (see {\em 3}) which might
make vertex $(i,*,j)$ isolated due to missing successors (see {\em 4}). 
Alternatively, vertex $(i,*,j)$ could become mergeable if removal of 
$((i,*,j),(j,*,k))$ makes its neighborhood identical with that of another 
vertex $(i,*).$ The two derivatives are added and one of the vertices is 
removed
from $\tilde{G}$ together with all incident edges. Similarly, vertex $(j,*,k)$
can become isolated or mergeable (see {\em 5}).
Figure~\ref{fig:bat}~(c) shows the result of eliminating 
$(0,2,3)$ from Figure~\ref{fig:bat}~(b).

A face elimination sequence is {\em complete} if the resulting line dag is tripartite
corresponding to a bipartite dag and representing the Jacobian of the
underlying program. 
Both termination and numerical correctness are shown in \cite{Naumann2004Oao}
as immediate consequences of the chain rule as stated in 
Equation~(\ref{eqn:cr}). Special cases can be derived including rules for
{\em edge} \cite{Naumann2004Oao} and {\em vertex elimination} 
\cite{Griewank1991OtC} in $G.$ For example, any complete face, edge or vertex elimination sequence transforms the dag in Figure~\ref{fig:bat}~(a) into a bipartite dag 
$G'=(V'_1 \cup V'_2,E'),$ where $V'_1=\{-1,0\},$ $V'_2=\{4,5,6\}$ and 
$E'=V'_1 \times V'_2.$

Assignment of a unit cost to single face eliminations yields 
the following combinatorial optimization problem.
\begin{Problem}[\sc Face Elimination]
Find the shortest complete face elimination sequence.
\end{Problem}
Obviously, computational cost should not be reduced to simple operations count
in general. Temporal and spatial locality properties of numerical programs
have a substantial effect. This fact will be accounted for in 
Section~\ref{sec:3}.

The restriction of the chain rule of differentiation to a structural 
formulation eliminates consideration of potential 
algebraic dependences (e.g., equality) amongst
entries of elemental Jacobians. Numerical values of the latter are not 
considered. The main argument of the proof of NP-hardness
of {\sc Jacobian Accumulation} is hence no longer applicable; see
\cite{Naumann2008OJa}. In fact, the
computational complexity of {\sc Face Elimination} remains an open problem.

Similar statements hold for the corresponding {\sc Edge} and {\sc Vertex 
Elimination} problems. 
The subgraph of Figure~\ref{fig:bat}~(a) induced by all vertices reachable from
$-1$ followed by re-indexing as $-1\rightarrow 0,$ $3\rightarrow 2,$ $4 \rightarrow 3,$ $5\rightarrow 4,$ and $6\rightarrow 5$
was shown in \cite{Naumann2002ETf} (and therein referred to as the {\em Lion}  
dag) to yield a lower computational cost (11 \fma) than the best vertex elimination 
sequence (12 \fma). Edge elimination was proposed together with corresponding greedy 
heuristics to solve this problem. Similar reasoning led to the definition
of face elimination driven by the observation that the 
cost of the best edge elimination sequence (23 \fma) can be improved for 
the {\em Bat} dag in Figure~\ref{fig:bat}~(a).
\begin{example} \label{ex:3}
The following face elimination sequence for the Bat dag induces a computational cost of $22$ \fma:
\begin{alignat*}{2}
F'(1,3,4) &\ass F'(3,4) \cdot F'(1,3) &&(1~\fma)\\
F'(2,3,6) &\ass F'(3,6) \cdot F'(2,3)  &&(1~\fma)\\
F'(-1,1,3,4) &\ass F'(1,3,4) \cdot F'(-1,1)  &\quad& (2~\fma)\\
F'(-1,1,3) &\ass F'(1,3) \cdot F'(-1,1)  &&(2~\fma)\\
F'(-1,1,3,5) &\ass F'(3,5) \cdot F'(-1,1,3)  &&(4~\fma)\\
F'(-1,1,3,6) &\ass F'(3,6) \cdot F'(-1,1,3)  &&(2~\fma)\\
F'(0,2,3,6) &\ass F'(2,3,6) \cdot F'(0,2)  && (2~\fma)\\
F'(0,2,3) &\ass F'(2,3) \cdot F'(0,2)  && (2~\fma)\\
F'(0,2,3,4) &\ass F'(3,4) \cdot F'(0,2,3)  && (2~\fma)\\
F'(0,2,3,5) &\ass F'(3,5) \cdot F'(0,2,3)  &&(4~\fma)\; .
\end{alignat*}
\end{example}
The scalar scenario has been studied for more than three decades.
Application of Equation.~(\ref{eqn:cr}) 
to graphs induced by single vertices (together with their
neighborhoods) yielded
vertex elimination which was suggested in \cite{Griewank1991OtC} 
alongside Markowitz-type 
greedy heuristics for determining a near-optimal elimination sequence.
Motivated by the Lion and Bat dags
edge \cite{Naumann2002ETf} and face 
\cite{Naumann2004Oao} elimination methods were developed. 
NP-hardness of the combinatorial {\sc Jacobian Accumulation} problem is 
shown in \cite{Naumann2008OJa}. Successful implementations of elimination 
techniques are described in \cite{Forth2004JCG,Tadjouddine2003HAD,Utke2008OAm}. 
Generalizations for higher derivatives were proposed \cite{Wang2016EPi}.
Contributions to the development of combinatorial optimization methods include 
\cite{Chen2012AIP,Lyons2012Rhf,Pryce2008FAD}.

\section{Generalized Face Elimination} \label{sec:3}

In real-world simulation scenarios
the Jacobians of nontrivial vector-valued multivariate elemental functions 
are typically not available.
{\em Elemental tangents}
\begin{equation} \label{eqn:T}
\begin{split}
\R^{n_j \times \dot{n}_j} \ni \dot{V}_j &\ass \dot{F}(i,j) \cdot \dot{V}_i
\equiv \dot{V}_j + F'(i,j) \cdot \dot{V}_i 
\end{split}
\end{equation}
and {\em elemental adjoints}
\begin{equation} \label{eqn:A}
\begin{split}
\R^{\bar{n}_i \times n_i} \ni \bar{V}_i &\ass  \bar{V}_j \cdot \bar{F}(i,j)
  \equiv       \bar{V}_i + \bar{V}_j \cdot F'(i,j)
\end{split}
\end{equation}
are assumed to be given instead, for example, as the result of the
application of an AD tool. 
W.l.o.g.,
we use vector mode (see Chapter 3 in \cite{Griewank2008EDP}) and incremental
notation for both tangents and adjoints allowing for a more compact typesetting
in the formulation of {\sc Rules}~\ref{r2} and \ref{r3}.
Vector tangent mode yields $F'(i,j)$ at the expense of $\mathcal{O}(n_i)$
evaluations of $F_j$ by setting $\dot{V}_i\ass I_{n_i}$ for $\dot{V}_j \ass 0.$
Vector adjoint mode yields $F'(i,j)$ at the expense of $\mathcal{O}(n_j)$
evaluations of $F_j$ by setting $\bar{V}_j\ass I_{n_j}$ for $\bar{V}_i \ass 0.$

Most AD tools generate scalar tangents 
$\dot{F}(i,j) \cdot \dot{\V}_i \in \R^{n_j}$ and adjoints
$\bar{\V}_i \cdot \bar{F}(i,j) \in \R^{1 \times n_i}.$ 
Hence and w.l.o.g., we expect realistic cost estimates for all scalar elemental tangents 
and adjoints to be provided by the user of the new method, e.g., obtained by 
profiling.

Moreover, tangents $\dot{F}_j \equiv F'_j \cdot (\dot{V}_i)_{i \prec j}$ and adjoints 
$\bar{F}_j \equiv \bar{V}_j \cdot F'_j$ associated with 
vertices are typically given instead of
$\dot{F}(i,j)$ and $\bar{F}(i,j)$ associated with individual edges in the dag. 
The latter can be obtained by specialization
of the former. Further technical details of this source code transformation
are omitted in favor of focusing on the combinatorial issues of AD mission
planning. For example, the AD tool Tapenade
\cite{Hascoet2013TTA} supports custom activity patterns for tangent and
adjoint subprograms. Similar functionality is provided by AD tools based
on operator and function overloading. See {\tt www.autodiff.org} for a 
collection of AD tools and a comprehensive bibliography on the general subject.

In this article, all complete face elimination sequences are assumed to have feasible memory
requirements. Conservatively, sufficient memory should be available for
storage of the inputs to all elemental functions in addition to the maximum
memory required for evaluating any single elemental adjoint. This
naive checkpointing \cite{Griewank1992ALG} scheme is likely to be feasible in practice given the limited size
of the dags under consideration. The additional evaluation of the primal should
be accounted for when defining the computational costs associated with the
vertices in the line dag.
A more general problem formulation
combining {\sc Face Elimination} and the computationally intractable 
{\sc DAG Reversal} \cite{Naumann2008DRi}
or {\sc Call Tree Reversal} \cite{Naumann2008CTR} problems
is the subject of ongoing research.

The search space of the combinatorial
{\sc generalized Face Elimination} problem to be formulated in the following
is induced by two elemental operations on the line dag: {\em Edge elimination}
can be applied to
edges connecting any two intermediate vertices, where at least one of them
has an elemental Jacobian; see {\sc Rule}~\ref{r2}. 
{\em Preaccumulation} is applicable to all intermediate vertices without 
elemental Jacobians; see {\sc Rule}~\ref{r3}.
\begin{Rule} \label{r2}
        {\em Edge elimination} is defined for all
        $((i,*,j),(j,*,k)) \in \tilde{E}_Z$
        with $F'(i,*,j) \neq \varnothing$ or
        $F'(j,*,k) \neq \varnothing$
        as follows:
\begin{enumerate}
        \item If $\nexists (i,*,k) \in \tilde{V} : $
\begin{itemize}
\item[] \hspace{-1mm} $P(i,*,k)=P(i,*,j)$ and $S(i,*,k)=S(j,*,k)$ 
\item[] \hspace{-1mm} or $F'(i,*,k)=\varnothing$ 
\end{itemize}
then {\em (fill)}
                \begin{itemize}
                        \item[] $\tilde{V} \ass \tilde{V} \cup (i,*,j,*,k)$
                        \item[] 
                  $\tilde{E} \ass \tilde{E} \cup~(P(i,*,j) \times (i,*,j,*,k)) $ 
\item[] $\qquad \quad \, \cup~((i,*,j,*,k) \times S(j,*,k)) $
                        \item[] $F'(i,*,j,*,k) \ass 0$ and $\dot{F}(i,*,j,*,k) \ass \varnothing$ and $\bar{F}(i,*,j,*,k) \ass \varnothing.$
                \end{itemize}
\item
        $F'(i,*,j,*,k) \pe \begin{cases} 
                \dot{F}(j,*,k)\,F'(i,*,j) $\ldots$ \\
\quad \text{if}~\dot{F}(j,*,k), F'(i,*,j) \neq \varnothing \\ 
                F'(j,*,k)\,\bar{F}(i,*,j) $\ldots$ \\
\quad \text{if}~F'(j,*,k), \bar{F}(i,*,j) \neq \varnothing \\
                F'(j,*,k) \cdot F'(i,*,j) $\ldots$ \\
\quad \text{if}~F'(j,*,k), F'(i,*,j) \neq \varnothing 
\end{cases}$
\item $\tilde{E} \ass \tilde{E} \setminus ((i,*,j),(j,*,k))$ {\em (remove edge)}.
\item If $S(i,*,j) = \varnothing$ then {\em (remove vertex)}
\begin{itemize}
        \item[] $\tilde{V} \ass \tilde{V} \setminus (i,*,j)$
        \item[] $\tilde{E} \ass \tilde{E} \setminus \left (P(i,*,j) \times (i,*,j) \right )$
\end{itemize}
                else
\begin{itemize}
\item[] if $\exists~(i,*) \neq (i,*,j) \in \tilde{V}: $
\begin{itemize}
\item[] $S(i,*) = S(i,*,j)$ and $F'(i,*) \neq \varnothing$
\end{itemize}
                then {\em (merge)}
\begin{itemize}
        \item[] $F'(i,*,j) \pe F'(i,*)$
\item[] $\tilde{V} \ass \tilde{V} \setminus (i,*)$
        \item[] $\tilde{E} \ass \tilde{E} \setminus ((i,*) \times S(i,*) \cup P(i,*) \times (i,*)).$
\end{itemize}
\end{itemize}
\item If $P(j,*,k) = \varnothing$ then {\em (remove vertex)}
\begin{itemize}
        \item[] $\tilde{V} \ass \tilde{V} \setminus (j,*,k)$
        \item[] $\tilde{E} \ass \tilde{E} \setminus \left ((j,*,k) \times S(j,*,k)\right )$
\end{itemize}
else
\begin{itemize}
\item[] if $\exists~(*,k) \neq (j,*,k) \in \tilde{V}:$
\begin{itemize}
\item[] $P(*,k) = P(j,*,k)$ and
                $F'(*,k) \neq \varnothing$
\end{itemize}
                then {\em (merge)}
\begin{itemize}
\item[] $F'(j,*,k) \pe F'(*,k)$
\item[] $\tilde{V} \ass \tilde{V} \setminus (*,k)$
\item[] $\tilde{E} \ass \tilde{E} \setminus ((*,k) \times S(*,k) \cup P(*,k) \times (*,k)).$
\end{itemize}
\end{itemize}
\end{enumerate}
\end{Rule}
{\sc Rule}~\ref{r2} extends {\sc Rule}~\ref{r1} with logic due to the elemental tangent-/adjoint-based
formulation of the chain rule. Fill-in yields new vertices with zero
Jacobians $F'(i,*,j,*,k)$ and without tangents nor adjoints (see {\em 1}).
The latter
are available for elemental functions represented by vertices in the
original line dag, exclusively.

The \fma\ operation comes in three flavors (see {\em 2}).
If both $F'(i,*,j)$ and $F'(j,*,k)$ are given, then their product\footnote{We use $\cdot$ to distinguish the matrix product from the Jacobian-free propagation of tangents or adjoints.} is added
to $F'(i,*,j,*,k).$ Otherwise, either the tangent of $(j,*,k)$ is seeded with
$F'(i,*,j)$ or the adjoint of $(i,*,j)$ with $F'(j,*,k).$ At least one of the
two elemental Jacobians is required to be given for $((i,*,j),(j,*,k))$
to qualify as eliminatable, which results in its
removal from $\tilde{G}$ (see {\em 3}). Merging (see {\em 4} and {\em 5}) is similar
to {\sc Rule}~\ref{r1} with the additional requirement that the elemental Jacobian of
the vertex to be merged needs to be available.

The cost 
of eliminating an eliminatable edge $((i,*,j),(j,*,k))$ 
according to the three alternatives in {\em 2} 
is equal to
$n_i \cdot \Ert(\dot{F}(j,*,k)),$ or
$n_k \cdot \Ert(\bar{F}(i,*,j))$ or 
$n_i n_j n_k,$ respectively. The cost of matrix multiplication is assumed to
be equal to the number of \fma~performed. The costs of matrix additions due to 
merging are neglected as it is dominated by the cost of {\em 2}.
%

\begin{Rule} \label{r3}
        {\em Preaccumulation}
        is defined for all $(i,j) \in \tilde{V}_Z$ with $F'(i,j) = \varnothing$ as follows:
\begin{enumerate}
        \item $F'(i,j)\ass \begin{cases}
                        \dot{F}(i,j)\,  I_{{n}_i} \\
                        I_{{n}_j} \, \bar{F}(i,j) 
        \end{cases}$
\item
        If $\exists~(i,*,j) \in \tilde{V}: $
\begin{itemize}
\item[] $P(i,*,j)=P(i,j)$ and $S(i,*,j)=S(i,j)$ and $F'(i,*,j) \neq \varnothing$
\end{itemize}
                then {\em (merge)}
\begin{itemize}
        \item[] $F'(i,j) \pe F'(i,*,j);$
        $\tilde{V} \ass \tilde{V} \setminus (i,*,j).$
\end{itemize}
\end{enumerate}
\end{Rule}
Preaccumulation is performed by seeding either the tangent or the adjoint 
with the corresponding identities (see {\em 1}). The less costly alternative
will be selected in the context of the combinatorial
{\sc Generalized Face Elimination} problem to be formulated further
below.
Mergeability of $(i,*,j)$ might follow under the same conditions as
in {\sc Rule}~\ref{r2} (see {\em 2}). 

Termination follows from an argument similar to that in \cite{Naumann2004Oao}. 
Correctness
of {\em 1} follows from the definition of tangents and adjoints. The merge rule
in {\em 2} preserves invariance under the chain rule.

The cost of preaccumulating $F'(i,j)$ is equal to 
$\min\{n_i \cdot \Ert(\dot{F}(i,j)),n_j \cdot \Ert(\bar{F}(i,j))\};$
the costs of matrix additions due to merging are neglected due to dominance
of the cost of {\em 1}.
In summary, the cost of a generalized face elimination amounts to the cost of the 
corresponding edge elimination added to the cost of potentially required 
preaccumulations; see also {\sc Example~\ref{ex:gfe}}.

A generalized version of {\sc Face Elimination} can be formulated
for given costs of all elemental tangents and adjoints.
\begin{Problem}[\sc Generalized Face Elimination]
Find a complete generalized face elimination sequence of minimal cost.
\end{Problem}
W.l.o.g., we assume unit \fma\ cost in matrix products. 
\begin{example} \label{ex:gfe}
The illustration of all aspects of {\sc Rules}~\ref{r2} and \ref{r3} with the 
help of examples turns out to be infeasible under the given space restrictions
due to combinatorial explosion. Instead, we consider the simplest possible 
scenario $F=F_2 \circ F_1$ with $n=n_1=10,$ $m_1=n_2=2,$ $m=m_2=5$, 
$\Ert(\dot{F}(0,1))=\Ert(\dot{F}(1,2))=1000$ and, optimistically,\footnote{The 
run time of the evaluation of an adjoint often exceeds that of a tangent by factors of more than two. The relative run time varies significantly for different 
AD tools.} 
$\Ert(\bar{F}(0,1))=\Ert(\bar{F}(1,2))=2000.$ 
\begin{alignat*}{2}
F'&={\dot{F}(1,2) \cdot (\dot{F}(0,1)  \cdot I_{n_1})} &&\Rightarrow \Ert=20.000 \\
F'&={\dot{F}(1,2) \cdot (I_{m_1} \cdot \bar{F}(0,1))} &&\Rightarrow \Ert=14.000 \\
F'&={(I_{m_2} \cdot \bar{F}(1,2)) \cdot \bar{F}(0,1)} &&\Rightarrow \Ert=20.000 \\
F'&={(\dot{F}(1,2) \cdot I_{n_2}) \cdot \bar{F}(0,1)} &&\Rightarrow \Ert=12.000 \\
F'&={(\dot{F}(1,2) \cdot I_{n_2}) \cdot (I_{m_1} \cdot \bar{F}(0,1))} &&\Rightarrow \Ert=6.100 \\
F'&={(I_{m_2} \cdot \bar{F}(1,2)) \cdot (I_{m_1} \cdot \bar{F}(0,1))} &&\Rightarrow \Ert=14.100 \\
F'&={(I_{m_2} \cdot \bar{F}(1,2)) \cdot (\dot{F}(0,1) \cdot I_{n_1})} &&\Rightarrow \Ert=20.100 \\
F'&={(\dot{F}(1,2) \cdot I_{n_2}) \cdot (\dot{F}(0,1) \cdot I_{n_1})} &&\Rightarrow \Ert=12.100 
\end{alignat*}

The optimal cost of $6.100$ is reached when preaccumulating $F'(0,1)$ in 
adjoint mode followed by preaccumulation of $F'(1,2)$ in tangent mode
and evaluation of the matrix product $F'(1,2) \cdot F'(0,1).$ It undercuts
the costs of the standard tangent (first line) and adjoint 
(third line) modes of AD by a factor of more than three.
\end{example}
{\sc Generalized Face Elimination} can be solved efficiently 
by dynamic programming
for single-path dags resulting from simple chains of elemental function 
evaluations \cite{Naumann2020OoG}. 
Neither sparsity \cite{Averick1994CLS} nor scarcity \cite{Griewank2005AaE} of
elemental Jacobians are taken into account by {\sc Rules}~\ref{r2} 
and \ref{r3}. 
Integration of both aspects will require additional annotation of the 
dag as well as modification of the elimination rules to detect and
keep track of sparsity and scarcity. Face elimination over sparse elemental 
Jacobians inherits its computational intractability from Jacobian Compression
through {\sc Graph Coloring} \cite{Gebremedhin2005WCI}.

{\sc Generalized Edge Elimination} (GEE) and 
{\sc Generalized Vertex Elimination} (GVE)
in dags follow naturally from
{\sc Generalized Face Elimination} (GFE). 
A hierarchy of search spaces is spanned such
that GVE $\subseteq$ GEE $\subseteq$ GFE. 

\section{Combinatorial Optimization}

The GFE problem yields a combinatorial search space of 
exponential size in the number of edges $|\tilde{E}|$ in the line dag 
$\tilde{G}$. Computational intractability can be proven for the practically
relevant scenario allowing for algebraic dependences amongst entries of
elemental Jacobians \cite{Naumann2008OJa}. The computational complexities of
the structural variants of GFE, GEE and GVE remain open. Ongoing research
puts focus on proving likely hardness and, hence, the development of 
appropriate methods including branch and bound and (greedy) heuristics. Reference implementations of all
methods described in the following can be found on
{\tt \href{https://github.com/STCE-at-RWTH/ADMission/tree/ACDA23}{github.com/STCE-at-RWTH/ADMission/tree/ACDA23}.}
All our numerical results can thus be reproduced.

\subsection{Heuristics}

The two obvious choices for generalized face elimination sequences 
correspond to sparse tangent \cite{Bischof1992AGD} 
and sparse adjoint \cite{Griewank2008EDP} modes of AD.
In sparse tangent mode, all elemental Jacobians associated with minimal vertices
are preaccumulated followed by the evaluation of the remaining elemental tangents in topological order.
Sparse adjoint mode preaccumulates elemental Jacobians of maximal vertices
followed by the evaluation of the remaining elemental adjoints in reverse 
topological order. For example, for the line dag in Figure~\ref{fig:bat}~(b)
sparse tangent mode would proceed as follows
\begin{align*}
F'(-1,1)&\pe \dot{F}(-1,1) \cdot I_{n_{-1}}\, ; \;\; F'(0,2) \pe \dot{F}(0,2) \cdot I_{n_0} \\
F'(-1,3)&\pe \dot{F}(1,3) \cdot F'(-1,1)\, ; \;\; \ldots
\end{align*}
starting from $F'(i,j)=0$ for $i=-1,\ldots 3$ and $j=1,\ldots 6.$ 

\noindent Sparse adjoint mode operates correspondingly:
\begin{align*}
F'(3,6)&\pe I_{m_6} \cdot \bar{F}(3,6)\, ; \;\; F'(2,6) \pe I_{m_6} \cdot \bar{F}(2,6)\; \ldots \\
F'(2,6)&\pe F'(3,6) \cdot \bar{F}(2,3)\, ; \; \ldots 
\end{align*}
Variants include locally optimal sparse tangent (similarly, adjoint) mode, where 
$F'(i,j)=F'(i,k) \cdot F'(k,j)$ allows preaccumulation of $F'(i,k)$ (in either
tangent or adjoint mode) followed
by the explicit matrix product in case the cumulative local cost undercuts
that of $F'(i,j)=\dot{F}(i,k) \cdot F'(k,j).$ 

Additionally, we have been experimenting with a range of  
heuristics. Results will be presented for a simple greedy heuristic based
on the local minimization of fill and accounting for 
negative fill due to merging. Ties are broken by either sparse tangent 
or adjoint modes.

\subsection{Branch and Bound}

A branch and bound algorithm was implemented as a deterministic method for
solving GFE. Its computational cost is proportional
to the size of the branch and bound search tree which
grows exponentially
with the size of the given (line) dag. Applicability is limited to
small dags. For example, the solution of GFE for the single Newton step
in Figure~\ref{fig:newton} takes less than a second on a standard workstation.
Two consecutive Newton steps yield a GFE instance whose branch and bound 
solution can take several weeks.
Partial exploration of infeasible search spaces becomes possible by letting 
the method run for a fixed time period, possibly combined with random 
multi-start.

We build on \cite{Mosenkis2018Olb}, where lower bounds for optimal 
Jacobian accumulation by vertex and edge elimination are derived. 
An upper bound on the cost of an optimal generalized face elimination sequence
is continuously updated by applying one of the greedy heuristics proposed 
above to the remainder graph. The resulting value is added to the cost of
the current partial generalized face elimination sequence corresponding
to the current position in the search tree. 

A lower bound
on the cost of subtrees of the search tree is obtained by adding lower 
bounds on the respective costs at which each of the remaining vertices can be 
eliminated from the line dag to the cost of the current partial generalized 
face elimination sequence. Highly conservative lower bounds are obtained by 
considering vertices of minimum 
dimensions which are connected to the current vertex by a directed path.
If the lower bound exceeds the current upper bound, then the current subtree
can be disregarded from further search.

Optimality-preserving compression of the search tree results from the 
exploitation of commutativity of mutually independent faces with respect to
computational cost. The elimination of
one face has no impact on the cost of elimination of the other in this case.
The same line dag results from the consecutive elimination of both faces 
in either order. Further compression results from
preaccumulation considered only in the context of enabling the elimination
of an incident edge in the line dag. Preservation of optimality follows 
immediately from the fact that the overall cost of a generalized face elimination sequence is invariant with respect to postponing preaccumulation until 
just prior to the elimination of the first incident edge.

The branch and bound method may have practical relevance in the context of 
AD mission planning, where a coarse-grain view of a numerical simulation 
code is likely to yield relatively 
small dags over expensive elemental functions. While previously introduced
heuristics can be expected to deliver satisfactory elimination sequences
in most cases,
certain 
high-end applications may benefit from additional savings in run time due to
a deterministically optimal solution.

\begin{example} \label{ex:4}
We continue {\sc Example~\ref{ex:2}}.
Tangents and adjoints are given for $F_1,F_2,$ and $F_3.$ Elapsed run times
of their scalar evaluation are attached to the respective edges in 
{\em Figure~\ref{fig:newton}}. 
Adjoints are twice as expensive as tangents. 
The cost vanishes for $(-1,3)$ as $F'(-1,3)=I_n$ is known. 
All elemental Jacobians are assumed to be dense except $F'(-1,3).$
For $n=m=10$
both (dense) tangent and adjoint modes yield a cost of
$124,000=20 \cdot (200+2,000+200+2,000+800+1,000) = 10 \cdot (400+4,000+400+4,000+1,600+2,000).$ 
The greedy heuristic yields
\begin{alignat*}{2}
F'(-1,2)&\ass\dot{F}(-1,2) \cdot I_{10} &\quad &\Rightarrow \Ert=2,000 \\
F'(-1,2,3)&\ass\dot{F}(2,3) \cdot F'(-1,2) &&\Rightarrow \Ert\pe 8,000 \\
F'(0,2)&\ass\dot{F}(0,2) \cdot I_{10} &&\Rightarrow \Ert\pe 2,000 \\
F'(0,2,3)&\ass\dot{F}(2,3) \cdot F'(0,2) &&\Rightarrow \Ert\pe 8,000 \\
F'(-1,1)&\ass\dot{F}(-1,1) \cdot I_{10} &&\Rightarrow \Ert\pe 20,000 \\
F'(-1,1,3)&\ass\dot{F}(1,3) \cdot F'(-1,1) &&\Rightarrow \Ert\pe 10,000 \\
F'(0,1)&\ass\dot{F}(0,1) \cdot I_{10} &&\Rightarrow \Ert\pe 20,000 \\
F'(0,1,3)&\ass\dot{F}(1,3) \cdot F'(0,1) &&\Rightarrow \Ert\pe 10,000 \\
F'(-1,3)&\pe F'(-1,1,3) &\;&\Rightarrow \Ert\pe 0 
\end{alignat*}
with a total cost of $80,000,$ which represents an improvement over 
dense AD by more than 35\%. Sparse adjoint mode turns out to be 
equivalent to dense adjoint mode due to the single output.

A global optimum of $74,000$ computed by branch and bound results in another 
minor improvement. The corresponding elimination sequence proceeds as follows:
\begin{alignat*}{2}
F'(0,2)&\ass \dot{F}(0,2) \cdot I_{10} &&\Rightarrow \Ert=2,000 \displaybreak[0] \\
F'(2,3)&\ass \dot{F}(2,3) \cdot I_{10} &&\Rightarrow \Ert \pe 8,000 \displaybreak[0] \\
F'(0,2,3)&\ass  F'(2,3) \cdot F'(0,2) &&\Rightarrow \Ert\pe 1,000 \displaybreak[0] \\
F'(-1,2)&\ass \dot{F}(-1,2) \cdot I_{10} &&\Rightarrow \Ert\pe2,000 \displaybreak[0] \\
F'(-1,2,3)&\ass F'(2,3) \cdot F'(-1,2) &\;&\Rightarrow \Ert\pe 1,000 \displaybreak[0] \\
F'(-1,1)&\ass \dot{F}(-1,1) \cdot I_{10} &&\Rightarrow \Ert\pe20,000 \displaybreak[0] \\
F'(-1,1,3)&\ass \dot{F}(1,3) \cdot F'(-1,1) &&\Rightarrow \Ert\pe 10,000 \displaybreak[0] \\
F'(0,1)&\ass \dot{F}(0,1) \cdot I_{10} &&\Rightarrow \Ert\pe20,000 \displaybreak[0] \\
F'(0,1,3)&\ass\dot{F}(1,3) \cdot F'(0,1) &&\Rightarrow \Ert\pe 10,000 \displaybreak[0] \\
F'(-1,3)&\pe F'(-1,1,3) &\;&\Rightarrow \Ert\pe 0 
\end{alignat*}
The total savings in cost add up to more than 40\% in comparison with
dense tangent and adjoint modes.
\end{example} 

\subsection{Validation}


An optimal bracketing of a generalized elemental Jacobian chain product can be 
computed
by dynamic programming; see \cite{Naumann2020OoG} for details.
We use the implementation described in this article for 
validation of the generalized face elimination rule on the resulting
single-path dags. 
Optimal elimination sequences computed by branch and bound 
reproduce the results of dynamic programming for all of a substantial number
of randomly generated test instances. 

Full numerical validation requires the (automated) implementation of the 
solutions of AD mission planning as actual derivative code. Subsequent
evaluation of Jacobians yields numerical values that can be validated, e.g.,
against finite difference approximations.
The development of the full AD mission pipeline is a major effort and
the subject of ongoing activities; see also Section~\ref{sec:concl}.

\section{Implementation}

A sample standard session of version 1.0 of the {\tt admission} software 
is presented.
More comprehensive information about building and using the AD mission planner can be
found in the public repository.
The software is currently developed with C++ for the Linux operating system
running on Intel-x86-based computers. 
Porting to other compute platforms is left as an option for future development.

The following three files play roles in AD mission planning sessions:
\begin{itemize}
\item[] \verb!admission!: the binary executable (name due to standard build process),
\item[] \verb!config.in!: the configuration file in text format (names may vary),
\item[] \verb!dag.xml!: the target dag in GraphML file format\footnote{http://graphml.graphdrawing.org/} (names may vary).
\end{itemize}
Here we assume all three of them to be stored in the current directory.
A session is run by typing
\begin{verbatim}
admission config.in
\end{verbatim}
on the command line. 

We take a closer look at the application of the branch and bound algorithm to 
the Lion dag. Its vertices and edges are defined 
in \verb!dag.xml!. For example, vertex $0$ becomes
\begin{verbatim}
<node id="0">
  <data key="index">0</data>
  <data key="vector_size">2</data>
</node> .
\end{verbatim}
Edge $(0,1)$ is described as
\begin{verbatim}
<edge id="0" source="0" target="1">
  <data key="tangent_cost">2</data>
  <data key="adjoint_cost">4</data>
  <data key="has_jacobian">0</data> ,
</edge>
\end{verbatim}
where the Boolean attribute \verb!has_jacobian! indicates whether the 
elemental Jacobian is available in addition to the elemental tangent and adjoint.

The contents of the configuration file
\begin{verbatim}
dag dag.xml
method BranchAndBound
\end{verbatim}
turn out to be self-explanatory. The following output is generated:
\begin{lstlisting}[basicstyle=\small,numbers=right,numberstyle=\tiny]
elapsed time: 0.00701579s

elimination sequence 
  (operation mode target cost):
  ACC TAN (0 1) 4
  ELI TAN (0 1 2) 2
  ELI TAN (0 1 4) 2
  ELI TAN (0 2 3) 6
  ELI TAN (0 2 4) 2

dense tangent cost: 16
dense adjoint cost: 64
optimized cost: 16

branch and bound statistics:
  number of nodes in search space: 416
  number of nodes visited: 190
\end{lstlisting}
An optimal elimination sequence preaccumulates the elemental Jacobian of edge $(0,1)$ in tangent mode (line 5) followed by tangent propagation through the 
remaining edges (lines 6-9). Both dense tangent and adjoint costs as well
as the optimum obtained by branch and bound (which happens to be equal to the 
cost of sparse tangent mode in this case) are displayed (lines 11-13). 
Additional statistics include the total number of nodes in the branch and 
bound tree (line 16) and the number of nodes actually visited (line 17).
Visits of the remaining nodes were avoided due to bounding. 

Further methods to be selected from are
sparse tangent (\verb!SparseTangent!), sparse adjoint
(\verb!SparseAdjoint!) and the greedy minimum fill (\verb!GreedyMinFill!) heuristic. Variants of the AD missions to be planned can be obtained by editing
the dag input file. For example, an instance with all elemental Jacobians available
can be generated by setting the \verb!has_jacobian! attribute of all edges to 
\verb!1!. Results for both classical (not generalized) versions of the 
Lion and Bat dags reported in Section~\ref{sec:2} were obtained this way.

\section{Results} \label{sec:cs}

We list results of all modes ordered as 
{\small
\verb!DenseTangent!, 
\verb!DenseAdjoint!, \verb!SparseTangent!, \verb!SparseAdjoint!, 
\verb!GreedyMinFill!, \verb!BranchAndBound!} 
for three case studies.

\subsection{Lion and Bat}

The following results are obtained
for the generalized versions of the Lion and Bat dags with 
respective edge costs as in Figure~\ref{fig:bat}~(a):
Lion: (16,64,16,46,20,16),
Bat: (48,96,32,64,36,32).
In both cases, sparse tangent mode turns out to be optimal.
Except in the case of branch and bound these results can easily be reproduced
by inspection of Figure~\ref{fig:bat}~(a).

\subsection{OpenFOAM}
One of our main target applications is OpenFOAM,
an open-source computational fluid dynamics solver suite, 
based on the finite volume method. We build on an algorithmic adjoint version
of OpenFOAM described in \cite{Towara2018Aam}.
The dag is obtained by instrumenting one step of the incompressible, 
steady state SIMPLE algorithm assuming laminar flow.
It is part of the public repository.
Velocities, face fluxes, and pressures are both in- and outputs to and from
the $i$-th iteration step of the SIMPLE algorithm.
The intermediate vertices correspond to objects allocated by OpenFOAM 
in order to solve the linearized Navier-Stokes equations.
Conceptually, the main tasks are the assembly and solution of the linearized 
momentum equations and the pressure correction equation.
A variety of dependencies exist between these two equations. An in-depth
description is beyond the scope of this article. Refer to the OpenFOAM
documentation\footnote{\tt www.openfoam.com/documentation} for further details.

Elemental tangent and adjoint costs were derived from the (very large) sizes 
(numbers of entries) of 
the {\em tapes} generated by the AD tool dco/c++ \cite{Naumann2016dDb}.
Elemental adjoints are approximately five times as expensive as elemental tangents.
The following results are obtained (rounded to the respective power of ten):
($8\cdot 10^9$, $41\cdot 10^9$, $5\cdot 10^9$, $27\cdot 10^9,$ $109\cdot 10^{12},$ $4.9\cdot 10^9$)
The greedy heuristic performs particularly poorly in this case. 
An improvement of 35\% over dense tangent mode was achieved by sparse tangent
mode. This solution was not improved significantly by branch and bound during
ten days using 24 threads on a local shared-memory computer. Only a very small 
fraction of the total search space 
($\approx 10^{-36}$) was explored at that point.


\subsection{Several Newton Steps}

Finally, we present the results for two additional variants of the Newton 
case study from {\sc Example~\ref{ex:Newton1}}. The evaluation of several Newton 
steps yields dags as concatenations of the dag of a single Newton step
from Figure~\ref{fig:newton}. For two steps we obtain
($248\cdot 10^3$, $248\cdot 10^3$, $200\cdot 10^3$, $248\cdot 10^3,$ $162\cdot 10^3,$ $152\cdot 10^3$). An improvement over dense AD
by nearly 40\% can be achieved by branch and bound running for 24 hours 
and using 48 threads on another local computer. 
Less than one percent of the search space was explored at that point.
The cost of an additional matrix product of
$2 \cdot 10^3=10 \cdot 10 \cdot 20$ needs to be added to the cost of 
preaccumulation
of two single Newton steps separately yielding a lower total cost of $150 \cdot 10^3=(2 \cdot 74 +2) \cdot 10^3.$ 
Similarly, ten steps yield
($1,240\cdot 10^3$, $1,240\cdot 10^3$, $2,600\cdot 10^3$, $1,240\cdot 10^3,$ $828\cdot 10^3,$ $812\cdot 10^3$) and hence an improvement of 35\% over standard 
AD. 
A tiny fraction of the total search space 
($\approx 10^{-127}$) was explored by branch and bound within 24 hours.
Again, separate optimal preaccumulation of single Newton steps yields a lower
cost of $(10 \cdot 74 + 9 \cdot 2) \cdot 10^3=758 \cdot 10^3.$
Ongoing investigations consider integration of this improved bound as part of a divide and conquer strategy within the branch and 
bound algorithm.

\section{Conclusion} \label{sec:concl}

Generalized face elimination is the basis of AD mission planning. Its formal
definition and illustration represents the main contribution of this article. 
Less emphasis was put on the development of corresponding combinatorial 
optimization methods. Promising first results obtained with relatively basic
variants of branch and bound and greedy heuristics motivate further ongoing and future activities in this area.

More generally, we envision a workflow combining (a) decomposition of 
a given numerical simulation into multivariate vector elementals, e.g. through 
annotation or instrumentation, (b) automatic derivation of the corresponding 
dag, (c) AD mission planning and (d) automatic generation of the optimized
derivative code. 
This article provides fundamental algorithmic techniques for the rigorous 
treatment of (c). 
Substantial further effort is required to address full
integration with (a), (b) and (d).

\end{document}